\documentclass{article}
 \usepackage[margin=1in]{geometry} 
 \usepackage{enumerate}
\usepackage{amsmath,amsthm,amssymb,amsfonts, fancyhdr, comment, parskip, tikz, caption, subcaption, graphicx}
\usepackage[colorlinks]{hyperref}
\usetikzlibrary{graphs}
\usetikzlibrary{shadings}
\usetikzlibrary{shapes.geometric}
\usetikzlibrary{arrows}
\usetikzlibrary{calc, positioning, fit, shapes.misc}

\definecolor{purple}{RGB}{148,0,211}
\definecolor{green}{RGB}{50,205,50}

\newtheorem{theorem}{Theorem}[section]

\newtheorem{lemma}[theorem]{Lemma}

\newcommand{\naturals}{\mathbb{N}}
\newcommand{\reals}{\mathbb{R}}
\newcommand{\TT}{{\cal T}}
\newcommand{\CC}{{\cal C}}

\newcommand{\GG}{{\cal G}}
\newcommand{\HH}{{\cal H}}

\newcommand{\PP}{{\cal P}}

\newcommand{\MM}{{\cal M}}

\newcommand{\EE}{{\cal E}}
\newcommand{\prob}{\mathbb{P}}
\newcommand{\eps}{{\varepsilon}}
\newcommand{\ex}{\text{ex}}
\newcommand{\E}{\mathbb{E}}
\newcommand{\usub}[2]{#1^{(#2)}}
\newcommand{\mc}[1]{\mathcal{#1}}
\newcommand{\rbrac}[1]{\left(#1\right)} 
\newcommand{\cbrac}[1]{\left\{ #1\right\}} 
\newcommand{\abrac}[1]{\left| #1\right|} 
\tikzstyle{arch} = [out=30, in=150]

\title{Generalized Ramsey numbers of cycles, paths, and hypergraphs}
\author{Deepak Bal\footnote{Department of Mathematics, Montclair State University, Montclair, NJ. \texttt{deepak.bal@montclair.edu}.}
\and Patrick Bennett\footnote{Department of Mathematics, Western Michigan University, Kalamazoo, MI. \texttt{patrick.bennett@wmich.edu}. The research of P. Bennett was partially supported by Simons Foundation Grant \#426894.}  
\and Emily Heath\footnote {Department of Mathematics, Iowa State University, Ames, IA. \texttt{eheath@iastate.edu}. The research of E. Heath was partially supported by NSF RTG Grant DMS-1839918.}
\and Shira Zerbib\footnote {Department of of Mathematics, Iowa State University, Ames, IA. \texttt{zerbib@iastate.edu}. The research of S. Zerbib was partially supported by NSF CAREER award no. 2336239,
NSF award no. DMS-195392,
and Simons Foundation award no. MP-TSM-00002629.}}
\date{}

\begin{document}

\maketitle

\begin{abstract}
    Given a $k$-uniform hypergraph $G$ and a set of $k$-uniform hypergraphs $\HH$, the generalized Ramsey number  $f(G,\HH,q)$ is   the minimum number of colors needed  to edge-color $G$ so that every copy of every hypergraph $H\in \HH$ in $G$ receives at least $q$ different colors.   In this note we obtain bounds, some asymptotically sharp, on several generalized Ramsey numbers, when $G=K_n$ or $G=K_{n,n}$ and $\HH$ is a set of cycles or paths, and when $G=K_n^k$ and $\HH$ contains a clique on $k+2$ vertices or a tight cycle. 
\end{abstract}

\section{Introduction}

In 1975, Erd\H{o}s and Shelah proposed the study of the following generalization of  Ramsey numbers: a \emph{$(p,q)$-coloring} of a graph $G$ is an edge-coloring of $G$ in which each copy of $K_p$ receives at least $q$ different colors. Let $f(n,p,q)$ be the minimum number of colors needed for a $(p,q)$-coloring of $K_n$.  Note that the special case $q=2$ is the classical multicolor Ramsey problem. These generalized Ramsey numbers were further developed by Erd\H{o}s and Gy\'arf\'as~\cite{EG}, and have been the focus of significant study over the years~\cite{axenovich2000,BEHK,BCDP,BED,CH1,CH2,CFLS,Mubayi1,mubayi2004}.  

In particular, Erd\H{o}s and Gy\'arf\'as gave an upper bound on $f(n,p,q)$ using the L\'ovasz Local Lemma, showing 
\[f(n,p,q)=O\left(n^{\frac{p-2}{\binom{p}{2}-q+1}}\right).\]
They proved that this bound was tight at the linear threshold $q=\binom{p}{2}-p+3$ and quadratic threshold $q=\binom{p}{2}-\left\lfloor\frac{p}{2}\right\rfloor+2$. However, Bennett, Delcourt, Li, and Postle~\cite{BDLP} recently showed that this bound can be improved by a logarithmic factor for all other values of $q$ (building on work of Bennett, English, and Dudek~\cite{BED} who proved it for about half of all possible values $q$). In fact, their improved upper bound applies in a more general setting: for subgraphs other than cliques, for hypergraphs of higher uniformity, and for list colorings. 

Given $k$-uniform hypergraphs $G$ and a family of $k$-uniform hypergraphs $\HH$, an  \emph{$(\HH,q)$-coloring} of $G$ is an edge-coloring of $G$ in which every copy in $G$ of every hypergraph $H\in \HH$ receives at least $q$ different colors. The minimum number of colors needed for an $(\HH,q)$-coloring of $G$ is denoted by $f(G,\HH,q)$. When $\HH$ contains only one hypergraph $H$, we write $(H,q)$-coloring and $f(G,H,q)$ instead of $(\HH,q)$-coloring and $f(G,\HH,q)$. Certain cases have received particular attention, including $G=K_{n,n}$ and $H=K_{p,p}$ (initiated by Axenovich, F\"uredi, and Mubayi~\cite{AFM}) and $G=K_n$ and $H=P_k$ (initiated by Krueger~\cite{K}). In~\cite{BDLP}, the authors showed that for a positive integer $k\geq2$, a $k$-uniform hypergraph $H$ with $|V(H)|>k$, a positive integer $q\leq |E(H)|$ such that $|E(H)|-q+1$ does not divide $|V(H)|-k$, and a $k$-uniform hypergraph $G$ on $n$ vertices, we have
\[f(G,H,q)=O\left(\left(\frac{n^{|V(H)|-k}}{\log n}\right)^{\frac{1}{|E(H)|-q+1}}\right),\] and moreover the same result holds in the setting of list colorings. 

This result was proved using a powerful new method in extremal combinatorics introduced independently by Delcourt and Postle~\cite{DP} under the name ``forbidden submatching method" and by Glock, Joos, Kim, K\"uhn, and Lichev~\cite{GJKKL} as the ``conflict-free hypergraph matching method."  Joos and Mubayi~\cite{JM} recently employed the latter version of this method to show that \[f(K_{n,n},C_4,3)=\frac{2}{3}n+o(n)\] and noted that a similar proof can be used in the non-bipartite case to show
\[f(K_n,C_4,3)=\frac{1}{2}n+o(n).\] In addition, they gave a shorter proof of a result of Bennett, Cushman, Dudek and Pra\l at~\cite{BCDP} showing that \[f(K_n, K_4, 5)=\frac{5}{6}n+o(n).\]

In this paper, we extend the results of~\cite{JM} in multiple ways. Our first result is to determine the asymptotics of the generalized Ramsey number for longer cycles in the complete graph at the linear threshold. In fact, we show that these same bounds hold when we forbid not just a single length of a cycle, but a family of cycles of consecutive lengths. For $4\leq k\leq \ell$, let $C_{[k,\ell]}=\{C_k,\ldots,C_{\ell}\}$. 

\begin{theorem}\label{thm:cycles}
    For fixed $k \ge 4$ and $\ell$ such that $k\leq \ell\leq \frac{\log \log n}{2 \log k}$, we have \[f\left(K_n,C_{[k,\ell]},3\right)=\frac{1}{k-2}n+o(n).\]
\end{theorem}

The lower bound in Theorem \ref{thm:cycles} follows only from the fact that in a $(C_k, 3)$-coloring we cannot have a monochromatic $P_k$, together with old results on $ex(n, P_k)$. 

Joos and Mubayi \cite{JM} showed that $f\left(K_{n,n},C_{4},3\right) = \frac 23 n + o(n)$. We generalize some of their work to give new upper bounds for $f\left(K_{n,n},C_{2k},3\right)$. (We also state a lower bound which follows similarly to the lower bound in Theorem \ref{thm:cycles}.)

\begin{theorem}\label{thm:bipartitecycles}
    For fixed $k\ge 3$, we have \[\frac{1}{2(k-1)}n +O(1) \leq f\left(K_{n,n},C_{2k},3\right)\leq\left(\frac{1}{2(k-1)} + \frac{1}{2a} \right)n+o(n),\]
    where 
    \begin{equation}
        a=a(k):= \left \lfloor \frac{2k-1+\sqrt{8k-7}}{2} \right \rfloor.
    \end{equation}
\end{theorem}
\noindent For large $k$ we have that $a$ is approximately $k$ and so the upper bound in Theorem \ref{thm:bipartitecycles} is approximately twice the lower bound. The proof of the upper bound in Theorem \ref{thm:bipartitecycles} is very similar to that in Theorem \ref{thm:cycles} so we just provide a sketch of the former. With a little more attention to detail we could also replace $C_{2k}$ in Theorem \ref{thm:bipartitecycles} with $C_{[2k, \ell]}$ for some $\ell = O(\log \log n)$ similarly to Theorem \ref{thm:cycles}. 

Next we generalize the result of Bennett, Cushman, Dudek and Pra\l at~\cite{BCDP} for $f(K_n,K_4,5)$ to the hypergraph setting. 

\begin{theorem}\label{thm:hypergraphs}
    We have \[f\left(K_n^k,K_{k+2}^k,\binom{k+2}{k}-1\right)=\frac{k^2+k-1}{k^2+k}n+o(n).\]
\end{theorem}

Finally, we give some results for paths in complete graphs, i.e. estimates for $f(K_n, P_k, q)$, where $P_k$ denotes the path on $k$ vertices.  Krueger \cite{K} proved that $f(K_n, P_k, q)$ is quadratic if and only if $q \ge \lceil k/2 \rceil +1$, and here we will focus on the quadratic threshold $q=\lceil k/2 \rceil +1$. It is not hard to see that $f(K_n,P_4,3)=f(K_n,P_5,4)=\binom{n}{2}$. We determine the asymptotics for the next two cases for paths at the quadratic threshold, and make some partial progress on the next case.

\begin{theorem}\label{thm:paths}
We have 
\begin{enumerate}
\item[(i)] $f(K_n,P_6,4)=\frac14 n^2+O(n)$,
\item[(ii)] $f(K_n,P_7,5)=\frac12 n^2+O(n)$, and 
\item[(iii)] $f_{\text{proper}}(K_n, P_8, 5) = \frac{7}{30}n^2 +O(n)$,
\end{enumerate}
where $f_{\text{proper}}(K_n, P_8, 5)$ is the minimum number of colors needed in a $(P_8, 5)$-coloring that is also a proper edge-coloring. 
\end{theorem}

In Section~\ref{sec:blackbox}, we present the main tools used in our proofs. In Section \ref{sec:cycles}, we prove Theorem \ref{thm:cycles}. In Section \ref{sec:cyclesbipartite}, we prove prove Theorem \ref{thm:bipartitecycles}. In Section \ref{sec:cliques}, we prove Theorem \ref{thm:hypergraphs}.  In Section \ref{sec:paths}, we prove Theorem \ref{thm:paths}.

\section{Tools}\label{sec:blackbox} 

In this section we introduce our two main tools for producing our colorings. 

\subsection{Lov\'asz Local Lemma}
We will make use of the following asymmetric version of the  Lov\'asz Local Lemma (see e.g., \cite{AS}).
Let $\mc{A}$ be a set of events and $G$ be a graph on vertex set $\mc{A}$. Denote by $N(A)$ the set of neighbors of $A$ in $G$. We say that $G$ is a  {\it dependency graph} for $\mc{A}$ if each event $A \in \mc{A}$ is mutually independent with its set of non-neighbors in $G$. 

\begin{lemma}[Lov\'asz Local Lemma] \label{lem:LLL}
Let $\mc{A}$ be a finite set of events in a probability space $\Omega$ and let $G$ be a dependency graph for $\mc{A}$. Suppose there is an assignment $x:\mc{A} \rightarrow [0, 1)$ of real numbers to $\mc{A}$ such that for all $A \in \mc{A}$ we have
\begin{equation}\label{eqn:LLLcond}
    \Pr(A) \le x(A) \prod_{B \in N(A)} (1-x(B)).
\end{equation}
Then, the probability that none of the events in $\mc{A}$ happen is 
\[
    \Pr\rbrac{\bigcap_{A \in \mc{A}} \overline{A}} \ge \prod_{A \in \mc{A}} (1-x(A)) >0.
\]
\end{lemma}

\subsection{Conflict-free hypergraph matchings}

Now we formally define the terminology necessary to state one of the main theorems of Glock, Joos, Kim, K\"uhn, and Lichev~\cite{GJKKL}. We opt for one of the most general versions of the theorem since it is necessary for some of our applications (e.g. to prove Theorem \ref{thm:cycles} with $\ell$ growing). 

Given a hypergraph $\HH$ and a vertex $v\in V(\HH)$, its \emph{degree} $\deg_{\HH}(v)$ 
is the number of edges in $\HH$ containing $v$. The maximum degree and minimum degree of $\HH$ are denoted by $\Delta(\HH)$ and $\delta(\HH)$, respectively. For $j\geq 2$, $\Delta_j(\HH)$ denotes the maximum number of edges in $\HH$ which contain a particular set of $j$ vertices, over all such sets. We use the notation $[\ell]:=\{1, \ldots, \ell\}$ and $[\ell]_m:=\{m, \ldots, \ell\}$.

In addition, for a (not necessarily uniform) hypergraph $\CC$ and an integer $k$, let $\CC^{(k)}$ be the set of edges in $\CC$ of size $k$. For a vertex $u\in V(\CC)$, let $\CC_u$ denote the hypergraph $\{C\backslash \{u\} \mid C\in E(\CC), u\in C\}$.

Given a hypergraph $\HH$, a hypergraph $\CC$ is a \emph{conflict system} for $\HH$ if $V(\CC)=E(\HH)$. A set of edges $E \subset \HH$ is \emph{$\CC$-free} if $E$ contains no subset $C\in \CC$. Given integers $d\geq 1$ and $\ell\geq 3$ and reals $\Gamma\geq 0$ and $\eps\in(0,1)$, we say $\CC$ is \emph{$(d,\ell,\Gamma, \eps)$-bounded} if  $\CC$ satisfies the following conditions:
\begin{enumerate}
    \item[(C1)]\label{cond:c1} $3\leq |C|\leq \ell$ for all $C\in\CC$;
    \item[(C2)]\label{cond:c2} $\sum_{j\in[\ell]} \frac{\Delta(\usub{\CC}{j})}{d^{j-1}}\leq \Gamma$ and~$| \{j\in[\ell]_3:\CC^{(j)}\neq\emptyset\} |\leq\Gamma$;
    \item[(C3)]\label{cond:c3} $\Delta_{j'}(\CC^{(j)})\leq d^{j-j'-\eps}$ for all $2\leq j'< j\leq \ell$. 
\end{enumerate}

Finally, given a $(d,\ell,\Gamma, \eps)$-bounded conflict system $\CC$ for a hypergraph $\HH$, we will define a type of weight function which can be used to guarantee that the almost-perfect matching given by Theorem~\ref{thm:blackbox} below satisfies certain quasirandom properties. We say a function $w:\binom{\HH}{j}\rightarrow[0,\ell]$ for $j\in\naturals$ is a \emph{test function} for $\HH$ if $w(E)=0$ whenever $E\in\binom{\HH}{j}$ is not a matching, and we say $w$ is \emph{$j$-uniform}. For a function $w:A\rightarrow \reals$ and a finite set $X\subset A$, let $w(X):=\sum_{x\in X} w(x)$. If $w$ is a $j$-uniform test function, then for each $E\subset \HH$, let $w(E)=w(\binom{E}{j})$. Given $j,d\in\naturals$, $\eps>0$, and a conflict system $\CC$ for hypergraph $\HH$, we say a $j$-uniform test function $w$ for $\HH$ is \emph{$(d,\eps,\CC)$-trackable} if $w$ satisfies the following conditions:
\begin{enumerate}
    \item[(W1)]\label{cond:w1} $w(\HH)\geq d^{j+\eps}$;
    \item[(W2)]\label{cond:w2} $w(\{{E\in \binom{\HH}{j}:E\supseteq E'\})\leq w(\HH})/d^{j'+\eps}$ for all $j'\in[j-1]$ and $E'\in\binom{\HH}{j'}$;
    \item[(W3)]\label{cond:w3} $|(\CC_e)^{(j')}\cap(\CC_f)^{(j')}|\leq d^{j'-\eps}$ for all $e,f\in \HH$ with $w(\{E\in\binom{\HH}{j}:e,f\in E\})>0$ and all $j'\in[\ell-1]$;
    \item[(W4)]\label{cond:w4} $w(E)=0$ for all $E\in\binom{\HH}{j}$ that are not $\CC$-free.
\end{enumerate}

\begin{theorem}[\cite{GJKKL}, Theorem 3.3]\label{thm:blackbox}
For all $r\geq 2$, there exists $\eps_0>0$ such that for all $\eps\in(0,\eps_0)$, there exists $d_0$ such that the following holds for all $d\geq d_0$.
		Suppose $\ell\geq 3$ is an integer and suppose $\Gamma\geq 1$ and $\mu\in(0,1/\ell]$ are reals such that $1/\mu^{\Gamma\ell}\leq d^{\eps^2}$.
		Suppose~$\HH$ is an $r$-uniform hypergraph on~$n\leq \exp(d^{\eps^2/\ell})$ vertices with~$(1-d^{-\eps})d\leq \delta(\HH)\leq\Delta(\HH)\leq d$ and~$\Delta_2(\HH)\leq d^{1-\eps}$ and suppose~$\CC$ is a~$(d,\ell,\Gamma,\eps)$-bounded conflict system for~$\HH$. Further, suppose~$\mc{W}$ is a set of~$(d,\eps,\CC)$-trackable test functions for~$\HH$ of uniformity at most~$\ell$ with~$\abrac{\mc{W}}\leq \exp(d^{\eps^2/\ell})$.
		Then, there exists a $\CC$-free matching~$\MM\subseteq \HH$ of size~$(1-\mu)n/r$ with~$w(\MM)=(1\pm d^{-\eps/900})\rbrac{|\MM|/|\HH|}^j w(\HH)$
		for all~$j$-uniform~$w\in\mc{W}$.
		
\end{theorem}

We will say that a hypergraph $\HH$ with $(1-d^{-\eps})d\leq \delta(\HH)\leq \Delta(\HH)\leq d$  is \emph{almost $d$-regular.}

\section{Cycles in complete graphs}\label{sec:cycles}

In this section we prove Theorem \ref{thm:cycles}. 

\subsection{Lower bound}\label{sec:cycleLB}

Suppose we have a $(C_k, 3)$-coloring of $K_n$. Then we cannot have any monochromatic $P_k$, and so each color class has at most $ex(n, P_k)$ edges. Faudree and Schelp \cite{FS} proved that for $n=q(k-1)+r$ and $0 \le r < k-1$ we have 
\begin{equation}\label{eqn:FS}
    ex(n, P_k)= q\binom{k-1}{2} + \binom{r}{2} = \frac{k-2}{2} n + o(n).
\end{equation}

Thus, the number of colors needed for a $(C_k, 3)$-coloring of $K_n$ is at least
\begin{equation}\label{eqn:cycleLB}
   \frac{\binom{n}{2}}{\frac{k-2}{2} n + o(n)} = \frac{n}{k-2} + o(n). 
\end{equation}

Thus we have the lower bound in Theorem \ref{thm:cycles}.

\subsection{Upper bound}\label{sec:cycleUB}

A very appealing aspect of the proof of Theorem \ref{thm:cycles} is how the proof of the lower bound illuminates what must be done to obtain a coloring for the upper bound. In particular, for the lower bound \eqref{eqn:cycleLB} to be asymptotically sharp, one would need a coloring in which almost all color classes have size asymptotically equal to the maximum possible number given in \eqref{eqn:FS}. The extremal construction for \eqref{eqn:FS} is given by $q$ vertex-disjoint copies of $K_{k-1}$ (and one copy of $K_r$ that is insignificant as $n$ grows). Thus we have the idea to color the graph using monochromatic copies of $K_{k-1}$ where two copies sharing a vertex cannot be the same color. 

We will construct a coloring of $K_n$ in two stages. The coloring in the first stage will use $\frac{n}{k-2}$ colors to color almost all the edges. To determine this coloring, we define appropriate hypergraphs $\HH$ and $\CC$ for which a $\CC$-free matching in $\HH$ corresponds to a partial coloring of $K_n$ with monochromatic copies of $K_{k-1}$. We will have to define our conflict system $\mc{C}$ so that each copy of $C_m$ for $m\in\{k,\ldots,\ell\}$ in the resulting coloring receives at least three colors. 

Theorem~\ref{thm:cyclesproperties} below summarizes the important properties of the coloring we will obtain from the first stage. In the second stage of the coloring, we then randomly color the remaining uncolored edges using a relatively small new set of colors. 

\begin{theorem}\label{thm:cyclesproperties}
For each $k$ there exists $\delta>0$ such that for all sufficiently large $n$ in terms of $\delta$, there is an edge-coloring of a subgraph $F\subset K_n$ with at most $\frac{n}{k-2}$ colors and the following properties:
\begin{enumerate}
    \item Every color class is a union of vertex-disjoint copies of $K_{k-1}$. \label{prop:cycle1}
    \item Every copy of $C_m$ in $F$ for $k\leq m\leq \ell$ receives at least three colors. \label{prop:cycle2}
    \item The graph $L=K_n\setminus E(F)$ has maximum degree at most $n^{1-\delta}$. \label{prop:cycle3}
    \item For each $xy\in E(K_n)$, the number of edges $x'y'\in E(L)$ such that $xy'$ and $yx'$ receive the same color is at most $n^{1-\delta}$. \label{prop:cycle4}
\end{enumerate}
\end{theorem}

\begin{proof}
We intend to apply Theorem \ref{thm:blackbox}, and we will use the following parameters: $r= \binom{k-1}{2}+(k-1)$, $\eps=\min\cbrac{\frac{\eps_0}{2}, \frac{1}{2k}}$, $\ell = \frac{\log \log n}{2 \log k}$,  $\Gamma = \ell k^\ell$, $\mu=1/\ell$, and $d=\binom{n}{k-2}$. Note that the condition $1/\mu^{\Gamma\ell}\leq d^{\eps^2}$ is satisfied since 
\[
\rbrac{\frac{1}{\mu}}^{\Gamma \ell} = \ell^{\ell^2 k^\ell} =\ell^{\ell^2 \sqrt{ \log  n}}= n^{o(1)}
\]
while 
\[
d^{\eps^2} = n^{(k-2)\eps^2 +o(1)}
\]
and $\eps$ depends only on $k$. In the future it will also be useful to remember that $d^{\eps} = n^{(k-2)\eps +o(1)}$ and that this power is $(k-2)\eps \le (k-2) \frac{1}{2k} < \frac12$.

Let $V=E(K_n)\cup\{v_i:v\in V(K_n),i\in [\frac{n}{k-2}]\}$. Let $\HH$ be the $r$-uniform hypergraph, where $r= \binom{k-1}{2}+(k-1)$, with vertex set $V$ which contains for each set $e$ of $k-1$ elements in $V(K_n)$ and each color $i\in [\frac{n}{k-2}]$ the edge $e_i:=\binom{e}{2}\cup\{v_i:v\in e\}.$ 

Note that a matching in $\HH$ corresponds to a set of edge-disjoint monochromatic copies of $K_{k-1}$ in $K_n$ in which no vertex appears in two cliques of the same color. Thus, a matching in $\HH$ yields an edge-coloring of a subgraph of $K_n$ with at most $\frac{n}{k-2}$ colors which satisfies property \ref{prop:cycle1} above. 

The hypergraph $\HH$ is almost $d$-regular with $d=\binom{n}{k-2}$. Indeed, for fixed $uv\in V$, we can count the edges in $\HH$ containing $uv$ by picking another $k-3$ vertices to form a clique with $uv$ in $\binom{n-2}{k-3}$ ways and picking the color of the clique in $\frac{n}{k-2}$ ways. For fixed $v_i\in V$, there are $\binom{n-1}{k-2}$ ways to pick an edge containing $v_i$ since we must pick the other $k-2$ vertices to complete the clique of color $i$ containing $v$. Thus, in both cases, the degree of the fixed vertex is $d-O(n^{k-3}) \ge \rbrac{1-d^{-\eps}}d$. 

Furthermore, we have $\Delta_2(\HH)\leq n^{k-3}< d^{1-\eps}$, since fixing a pair of vertices in $\HH$ fixes either at least three vertices of $K_n$ or at least two vertices of $K_n$ and the color of a monochromatic clique. 

 Note that property \ref{prop:clique1} is guaranteed by the fact that we are finding a matching in $\HH$ and the edges in a matching are disjoint. When one encodes a problem as a conflict-free hypergraph matching, there are fundamentally two ways to enforce that the resulting structure avoids violating certain constraints: one can encode those constraints so that they would correspond to intersecting edges of $\HH$ (i.e. violations of the matching property), or one can encode these constraints as conflicts. It is often beneficial to encode as many constraints as possible so that they are enforced by the matching property, since we would like our conflict system to be as sparse as possible so it satisfies the conditions of Theorem \ref{thm:blackbox}. Next we define a conflict system $\CC$ for $\HH$ with conflicts of even size $s$ for $4 \le s \le \ell$.  Let $V(\CC)=E(\HH)$. 
The conflicts in $\CC$ arise from 2-colored cycles $C$ in $K_n$ of length $m\in \{k,\ldots, \ell\}$. Note that any such 2-colored cycle must consist of an even number of alternately colored monochromatic paths, and each one of these paths must be contained in the vertex set of one of the monochromatic cliques in our matching in $\HH$. Thus a conflict will be a set of edges in $\HH$ which correspond to a collection of cliques in $K_n$ with a cyclic ordering such that the cliques alternate between two colors in this ordering, and such that each clique shares exactly one vertex each with the cliques before and after it. If we include these conflicts in our system, then a conflict-free matching in $\HH$ will indeed correspond to a partial edge-coloring of $K_n$ with none of the bicolored cycles we need to avoid. So we let $\CC$ be the set of such conflicts. 


Now we address Condition (C2). To this end, fix $e_{i}\in V(\CC)$. Since our conflicts are linear hypergraphs in $\HH$, we can determine a conflict of size $s$ containing $e_i$ by choosing, for each of the first $s-2$ edges (other than $e_i$),  $k-2$ new vertices together with one of the $k-1$ vertices from the last edge  to complete a clique in $\HH$ and then for the last edge picking just $k-3$ new vertices together with one vertex (of $k-1$) from each of the two edges it intersects. In addition, we must pick the second color $j\in[\frac{n}{k-2}]$ for the conflict. Thus, $\Delta(\CC^{(s)})\leq (k-1)^{s-1}\binom{n}{k-2}^{s-2}\cdot \binom{n}{k-3}\cdot \frac{n}{k-2}\le  k^\ell d^{s-1}$. Thus we have
\[
\sum_{s\in[\ell]} \frac{\Delta(\usub{\CC}{s})}{d^{s-1}} \le  \ell k^\ell = \Gamma.
\]
We also trivially have $|\{j\in[\ell]_3:\CC^{(j)}\neq\emptyset\} |\le \ell < \Gamma.$ Thus (C2) is satisfied.

We move on to checking Condition (C3). If we fix two edges, then their codegree is maximized when the edges either intersect in a vertex or share a color. In the first case, there are $O\left((n^{k-2})^{s-3}\cdot n^{k-3}\right)$ conflicts containing both edges, and in the second, there are $O\left((n^{k-2})^{s-4}\cdot(n^{k-3})^2\cdot n\right)$ such conflicts. Hence, $\Delta_2(\CC^{(2)})=O(d^{s-2-\eps})$. Similarly, the codegree of a fixed set $S$ of $s'\geq 3$ edges in $\HH$ is maximized when all $s'$ edges form a linear path in $\HH$, in which case there are at most $O((n^{k-2})^{s-s'-1}\cdot n^{k-3})$ conflicts containing $S$. Thus, we have $\Delta_{s'}(\CC^{(s)})=O(d^{s-s'-\eps})$. So, $\CC$ is a $(d,\ell, \Gamma, \eps)$-bounded conflict system for $\HH$.

Applying Theorem~\ref{thm:blackbox} with $\HH$ and $\CC$ would yield a conflict-free matching $M$ in $\HH$ which would correspond to a coloring of a subgraph $F$ of $K_n$ satisfying properties \ref{prop:cycle1} and \ref{prop:cycle2} of Theorem~\ref{thm:cyclesproperties}. 
In order to obtain a coloring which also satisfies properties \ref{prop:cycle3} and \ref{prop:cycle4} of Theorem~\ref{thm:cyclesproperties}, we now introduce appropriate test functions.

First, we define a 1-uniform test function for each vertex $v\in V(K_n)$. Let $S_v$ be the set of edges incident to $v$ in $K_n$. Let $w_v:E(\HH)\rightarrow\{0,k-2\}$ be the weight function which assigns to each edge $e_i\in E(\HH)$ the size of its intersection with $S_v$. Then we have  
\[w_v(\HH)=\sum_{uv\in S_v}\deg_{\HH}(uv)=nd-O(n^{k-3})\geq d^{1+\eps},\]
satisfying condition (W1). In addition, conditions (W2)-(W4) are trivially satisfied since $w_v$ is 1-uniform, so $w_v$ is a $(d,\eps,\CC)$-trackable  test function for all $v\in V(K_n)$. 

Now fix $\delta<\eps^3\log_n(d) = (1+o(1))(k-2)\eps^3$. Applying Theorem~\ref{thm:blackbox} with these test functions yields a $\CC$-free matching $M$ in $\HH$ for which 
\[w_v(M)>(1-d^{-\eps^3})d^{-1}w_v(\HH)>(1-n^{-\delta})n\]
for each $v\in V(K_n)$. 
Consequently, any $v\in V(K_n)$ has degree at most $n-(1-n^{-\delta})n=n^{1-\delta}$ in the uncolored subgraph $L\subset K_n$, proving property \ref{prop:cycle3} of Theorem~\ref{thm:cyclesproperties}.

In order to also obtain property \ref{prop:cycle4} of Theorem~\ref{thm:cyclesproperties}, we will define two more test functions for each edge $xy\in E(K_n)$ that will allow us to bound the number of edges $x'y'\in E(L)$ such that $xy'$ and $x'y$ receive the same color in $F$. For the first, let 
\[\PP_{xy}=\left\{\{e_i,e_i'\}:i\in\left[\frac{n}{k-2}\right], V(e)\cap V(e')=\emptyset,
x\in V(e), y\in V(e')\right\}.\] Let $w_{xy}$ be the indicator weight function for this set. This is a 2-uniform test function for $\HH$, and we will show that $w_{xy}$ is $(d,\eps,\CC)$-trackable. Condition (W1) holds for $w_{xy}$ since
\[w_{xy}(\HH)=|\PP_{xy}|\geq \binom{n}{k-2}^2\cdot\frac{n}{k-2}>d^{2+\eps}.\] 
In addition, for a fixed $\HH$-edge $e_i$, there are at most $O(n^{k-2})<\frac{w_{xy}(\HH)}{d^{1+\eps}}$ pairs in $\PP_{xy}$ which contain $e_i$, so condition (W2) is satisfied. Finally, for any two $\HH$-edges $e_i$ and $e_i'$ which form a pair in $\PP_{xy}$, there are at most $O(k^\ell(n^{k-2})^{s-2}n^{k-3}\cdot \frac{n}{k-2}\cdot n^{-1})=O(k^\ell n^{(s-1)(k-2)-1})<d^{s-1+\eps}$ sets of $s-1$ $\HH$-edges which complete a conflict with $e_i$ and, separately, with $e_i'$, for each $s\in\{4,6,\ldots,\ell\}$. This shows that condition (W3) holds. Finally, condition (W4) is vacuously true. Thus, $w_{xy}$ is $(d,\eps,\CC)$-trackable for each $xy\in E(K_n)$. 

For each $xy\in E(K_n)$, we also define a second test function $w'_{xy}$ which is an indicator weight function for the set 
\[\TT_{xy}=\left\{\{e_i,e_i',f_j\}:\{e_i,e'_i\}\in\PP_{xy}, j\in\left[\frac{n}{k-2}\right]\backslash\{i\}, |V(f)\cap V(e)|=1, |V(f)\cap V(e')|=1\right\}.\] 
This is a 3-uniform test function, which we now show is also $(d,\eps,\CC)$-trackable. Note that 
\[w'_{xy}(\HH)=|\TT_{xy}|=(k-2)^2(d\pm O(n^{k-3}))|\PP_{xy}|>d^{3+\eps},\]
so condition (W1) is met for $w_{xy}'$. Now fix an $\HH$-edge $e_i$. There are at most $O(n^{k-2}\cdot n^{k-3}\cdot \frac{n}{k-2})<\frac{w'_{xy}(\HH)}{d^{1+\eps}}$ triples in $\TT_{xy}$ containing $e_i$. Furthermore, if we fix two $\HH$-edges $e_i$ and $f_j$ of distinct colors, then there are at most $O(n^{k-3})<\frac{w'_{xy}(\HH)}{d^{2+\eps}}$ triples in $\TT_{xy}$ containing them both, and if we fix two $\HH$-edges $e_i$ and $e_i'$ of the same color, then there are at most $O(n^{k-4}\cdot \frac{n}{k-2})$ triples containing them both. Thus, condition (W2) is also satisfied for $w'_{xy}$. To check condition (W3), fix two $\HH$-edges which are in some triple together in $\TT_{xy}$. Either they share the same color or they share a vertex. In both cases, there are at most $O(k^\ell (n^{k-2})^{s-2}n^{k-3}\cdot \frac{n}{k-2}\cdot n^{-1})=O(k^\ell n^{(s-1)(k-2)-1})<d^{s-1+\eps}$ sets of $s-1$ $\HH$-edges which complete a conflict with either edge for each $s\in\{4,6,\ldots,\ell\}$. Finally, as before, condition (W4) is vacuously true, and $w'_{xy}$ is $(d,\eps,\CC)$-trackable.

We can now apply Theorem~\ref{thm:blackbox} with all of the test functions of type $w_v$, $w_{xy}$, and $w'_{xy}$ to obtain a $\CC$-free matching $M$ in $\HH$ for which property \ref{prop:cycle3} of Theorem~\ref{thm:cyclesproperties} holds and we have 
\[w_{xy}(M)=\left|\binom{M}{2}\cap \PP_{xy}\right|\leq (1+d^{-\eps^3})d^{-2}w_{xy}(\HH)\leq (1+n^{-\delta})\frac{n}{k-2}\] 
and 
\[w'_{xy}(M)=\left|\binom{M}{3}\cap \TT_{xy}\right|\geq (1-d^{-\eps^3})d^{-3}w'_{xy}(\HH)\geq (1-n^{-\delta})(k-2)n\]
for all $xy\in E(K_n)$.
Observe that for $xy\in E(K_n)$, each $x'y'\in E(L)$ with endpoints in distinct $\HH$-edges of the same color corresponds to an $\HH$-edge in a triple of $\TT_{xy}$ extending some pair from $\PP_{xy}$. Thus, the number of edges $x'y'\in E(L)$ as described in property \ref{prop:cycle4} of Theorem~\ref{thm:cyclesproperties} is at most the difference 
\[(k-2)^2 w_{xy}(M)-w'_{xy}(M)\leq n^{1-\delta},\]
which completes the proof of Theorem~\ref{thm:cyclesproperties}.
\end{proof}

To complete the proof of Theorem \ref{thm:cycles}, we now describe the second stage of our coloring. Assume we have completed the first stage and that Theorem~\ref{thm:cyclesproperties} holds with $2\delta$ in place of $\delta$, so we obtain a coloring with $\frac{n}{k-2}$ colors of a subgraph $F\subset K_n$ with the four desired properties. For the second stage, we color the remaining edges $L=E(K_n)\backslash F$ randomly with a new set $P$ of $c=n^{1-\delta}$ colors so that each edge in $L$ receives a color from $P$ with probability $1/c$, independently of the other edges. We will show using the Local Lemma that there is a coloring of $L$ which, when combined with the coloring of $F$, has no 2-colored cycle $C_m$ for any $m\in \{k,\dots,\ell\}$. 

We define three types of bad events. First, for any pair of adjacent edges $e,f\in L$, let $A_{e,f}$ be the event that both edges $e$ and $f$ receive the same color. Note that $\prob[A_{e,f}]=c^{-1}$. Next, for any cycle $D$ in $L$ of length $m\in\{k,\ldots,\ell\}$, let $B_D$ be the event that $D$ is properly-colored with 2 colors from $P$. Then $\prob(B_D)\leq c^{-m+2}$. Last, for each cycle $D$ in $K_n$ of length $m\in\{k,\ldots,\ell\}$ with $t\ge 1$ edges from $L$, let $C_{D}$ be the event that $D$ is 2-colored (this can only happen if all the edges from $L$ get the same color). Then $\prob[C_{D,i}]\leq c^{-t+1}$.  

Let $\EE$ be the collection of all bad events of each of these three types. Two events are mutually independent if they are edge-disjoint, meaning that they have no edges in $L$ in common. Fix an event $E\in\EE$. There are at most $\ell$ ways to pick a graph edge $e$ in the event $E$, and we will bound the number of events of each type which contain $e$. By property \ref{prop:cycle3} of Theorem~\ref{thm:cyclesproperties}, the number of events of type $A_{e,f}$ which contain $e$ is at most $2 \Delta(L)\leq 2n^{1-2\delta}$. For $3 \le m \le \ell$ there are at most $(\Delta(L))^{m-2}=n^{(m-2)(1-2\delta)}$ events of type $B_D$ which contain $e$ involving a cycle of length $m$. Finally, for each $2 \le t \le \ell$ there are at most $(k-2)^{t-1}(\Delta(L))^{t-1}\le k^\ell n^{(t-1)(1-2\delta)} $ events of type $C_{D}$ which contain $e$ involving $t$ edges from $L$.

We apply Lemma \ref{lem:LLL}. To each event $A_{e, f}$ we assign $x_A= 2c^{-1}= 2n^{-1+\delta}$, to each event $B_D$ involving a cycle of length $m$ we assign $x_{B, m} = 2c^{-m+2}= 2n^{-(m-2)(1-\delta)}$, and to each event $C_D$ involving $t$ edges from $L$ we assign $x_{C, t} = 2 c^{-t+1} = 2n^{-(t-1)(1-\delta)}$. Let $E$ be any particular bad event. Then 
\begin{align*}
    \prod_{B \in N(E)} (1-x(B)) & \ge 1 - \sum_{B \in N(E)} x(B)\\
    & \ge 1 - \ell \left[2n^{1-2\delta} \cdot x_A + \sum_{m=3}^{\ell} n^{(m-2)(1-2\delta)} x_{B, m} + \sum_{t=2}^{\ell} k^\ell n^{(t-1)(1-2\delta)} x_{C,t}\right]\\
    & = 1 - \ell \left[ 4n^{-\delta} + 2\sum_{m=3}^{\ell} n^{-(m-2)\delta}  + 2\sum_{t=2}^{\ell} k^\ell n^{-(t-1)\delta}\right]\\
    & = 1-o(1).
\end{align*}

Thus, Lemma \ref{lem:LLL} tells us that with positive probability none of the bad events happen, i.e. there exists a suitable coloring for the second stage. This completes the proof of Theorem \ref{thm:cycles}.

\section{Cycles in complete bipartite graphs}\label{sec:cyclesbipartite}

In this section we prove Theorem \ref{thm:bipartitecycles}. The lower bound follows easily from Gy\'{a}rf\'{a}s, Rousseau and Schelp \cite{GRS}, who proved the following. Suppose $G \subseteq K_{n, n}$ has no path on $2k$ vertices, where $k$ is fixed and $n \rightarrow \infty$. Then $G$ has at most $2(k-1)n + O(1)$ edges. Actually they proved a much more precise, non-asymptotic result but we will not use it. 

From here, the lower bound in Theorem \ref{thm:bipartitecycles} follows easily. Suppose we color the edges of $K_{n, n}$ such that every $C_{2k}$ has at least 3 colors. Then there is no monochromatic $P_{2k}$, and so each color class has at most $2(k-1)n + O(1)$ edges. Thus the number of colors is at least $n^2 / [2(k-1)n + O(1)]$, completing the proof of the lower bound.

\subsection{Upper bound}

Similarly to the upper bound in Theorem \ref{thm:cycles}, we will construct our coloring in two stages. In the first stage we will use a set $C$ of $\left(\frac{1}{2a} + \frac{1}{2(k-1)}\right)n$ colors to color almost all of the edges of $K_{n, n}$, and in the second stage we will use $o(n)$ colors for the rest of the edges of $K_{n, n}$. Say the bipartition of $K_{n, n}$ is $A \cup B$. 

In the coloring for the first stage, each color class will consist of vertex-disjoint copies of $K_{a, k-1}$, and furthermore these copies of $K_{a, k-1}$ will always share at most one vertex even if they are different colors. We will use Theorem \ref{thm:blackbox} to obtain this coloring, and to show that it has nice properties that will help us obtain the coloring in the second stage. 

A key difference from the proof of Theorem~\ref{thm:cycles} is the use of randomness in our construction of the hypergraph $\HH$. We will use the following concentration inequality due to McDiarmid~\cite{mcdiarmid} in our proof. 
\begin{theorem}\label{thm:McDiarmid}
Suppose $X_1,\ldots,X_m$ are independent random variables. Suppose $X$ is a real-valued random variable determined by $X_1,\ldots, X_m$ such that changing the outcome of $X_i$ changes $X$ by at most $b_i$ for all $i\in [m]$. Then, for all $t\geq 0$, we have \[\prob[|X-\E[X]|\geq t]\leq 2\exp\left(-\frac{2t^2}{\sum_{i\in [m]}b_i^2}\right).\]
\end{theorem}

Set 
\[
p:= \frac{a-1}{2(k-1)} + \frac{k-2}{2a}.
\]
By our choice of $a = \left \lfloor \frac{2k-1+\sqrt{8k-7}}{2} \right \rfloor$, we have $p \le 1$. Now we make a random set $E'$ of nonedge pairs $\{u, v\} \notin E(K_{n, n})$ by including each possible pair independently with probability $p$. When we pack vertex-disjoint copies of $K_{a, k-1}$, the nonedges of each copy of $K_{a, k-1}$ we include will always be in $E'$. Now set $V=E(K_{n, n})\cup E' \cup \{v_i:v\in A \cup B ,i\in C\}$. Let $\HH$ be the hypergraph on $V$ with the following edges. For each copy $K$ of $K_{a, k-1}$ in $K_{n, n}$ such that $\binom{V(K)}{2} \setminus E(K_{n, n}) \subseteq E'$  and each color $i\in C$, $\HH$ has the edge $\binom{V(K)}{2}\cup\{v_i:v\in V(K)\}.$ 

Then with high probability (with respect to our random choice of $E'$), $\HH$ is almost $d$-regular with 
\[
d=2\binom{n}{a-1} \binom{n}{k-2}|C|p^{a+k-1}.
\]
Indeed, for fixed $uv\in E(K_{n, n})$, we can count the edges in $\HH$ containing $uv$ by first choosing one of two possible orientations for our copy of $K_{a, k-1}$, then choosing another $a-1$ vertices from one side of the bipartition and $k-2$ from the other side, choosing a color and finally multiplying by the probability that all the required pairs are in $E'$. (Concentration easily follows from McDiarmid's inequality.) Now if instead we have $uv \in E'$, then we can complete an edge of $\HH$ as follows: depending on the orientation of our copy of $K_{a, k-1}$, we have either $\binom{n-2}{a-2} \binom{n}{k-1}$ or $\binom{n}{a} \binom{n-2}{k-3}$ choices for the rest of its vertices. We also choose a color and multiply by only $p^{a+k-2}$ since we are given that the pair $uv \in E'$, for a total of 
\[
\left(\binom{n-2}{a-2} \binom{n}{k-1}+\binom{n}{a} \binom{n-2}{k-3} \right)|C|p^{a+k-2} = (1+o(1))d
\]
by our choice of $p$. The final case to check for the almost $d$-regularity of $\HH$ is for a vertex of the form $v_i$. The degree of such a vertex is 
\[
\left(\binom{n-1}{a-1} \binom{n}{k-1}+\binom{n}{a} \binom{n-1}{k-2} \right)p^{a+k-1}= (1+o(1))d
\]
by our choice of $|C|$. 

Furthermore, we have $\Delta_2(\HH)\leq n^{a+k-3}$, since fixing a pair of vertices in $\HH$ fixes either at least three vertices of $K_{n, n}$ or at least two vertices of $K_{n, n}$ and a color. 

Next we define a conflict system $\CC$ for $\HH$ with edges of size $2s$ for $2 \le s \le k$.  Let $V(\CC)=E(\HH)$, and let the edges of $\CC$ correspond to copies of $C_{2k}$ in $K_{n, n}$ formed by two monochromatic matchings of size $s$. Now we address Condition (C2). To this end, fix an edge $e$ of $\HH$. Since our conflicts are linear hypergraphs in $\HH$, we can complete a conflict of size $s$ by choosing, for each of the first $s-2$ edges, $a+k-2$ new vertices together with one of the $a+k-1$ vertices from the last edge to complete a clique in $\HH$ and then for the last edge picking just $a+k-3$ new vertices together with one vertex (of $a+k-1$) from each of the two edges it intersects. In addition, we must pick the second color $j\in C$ for the conflict. Thus, $\Delta(\CC^{(s)})= O\left(\binom{n}{a+k-2}^{s-2}\cdot \binom{n}{a+k-3}\cdot |C| \right) =  O(d^{s-1})$. Thus for some $\Gamma = O(1)$ we have
\[
\sum_{s \le k} \frac{\Delta(\usub{\CC}{s})}{d^{s-1}} \le \Gamma.
\]
We also trivially have $|\{j:\CC^{(j)}\neq\emptyset\} | = O(1).$ Thus (C2) is satisfied. We omit the details of checking (C3) except to note that we will have $\Delta_{j'}(\CC^{(j)}) = O(d^{j-j'}/n)$ in all relevant cases. Thus we are able to apply Theorem \ref{thm:blackbox}.

The remainder of the proof is an application of the Local Lemma, for which one needs to verify certain properties of the coloring obtained from Theorem \ref{thm:blackbox} by defining appropriate test functions. Since this is very similar to the end of the proof of Theorem \ref{thm:cycles} we only give a brief sketch. By adapting the test functions from Theorem \ref{thm:cycles} we can establish that Theorem \ref{thm:blackbox} outputs a coloring where for each vertex $v$ in $K_{n, n}$, the number of uncolored edges $e$ of $K_{n, n}$ incident with $v$ is at most $n^{1-\delta}$ for some appropriately chosen $\delta>0$. Likewise for each uncolored edge $xy$ of $K_{n, n}$ there are at most $n^{1-\delta}$ other uncolored edges $x'y'$ of $K_{n, n}$ such that our coloring from Theorem \ref{thm:blackbox} contains two monochromatic copies of $K_{a, k-1}$ of the same color, one containing $xx'$ and the other containing $yy'$. With these properties in hand, we can easily adapt the application of the Local Lemma from Theorem \ref{thm:cycles} to our current setting, completing the proof of Theorem \ref{thm:bipartitecycles}. 

\section{Cliques in complete hypergraphs}\label{sec:cliques}

In this section we prove Theorem \ref{thm:hypergraphs}. 

\subsection{Lower bound}
Consider a $\left(K_{k+2}^k,\binom{k+2}{k}-1\right)$-coloring of $K_n^k$. Let $x_0$ be the number of pairs $(S, c)$ where $S$ is a set of $k-1$ vertices which is not contained in any edge of color $c$. Let $x_1$ be the number of pairs $(e, c)$ where $e$ is a $c$-colored edge such that no other $c$-colored edge intersects $e$ in $k-1$ vertices. Let $x_2$ be the number of pairs $(e_1 \cup e_2, c)$ such that $e_1$ and $e_2$ are $c$-colored edges which share $k-1$ vertices.  Then the total number of edges in $K_n^k$ is
\begin{equation}\label{eqn:edgecount}
x_1+2x_2 = \binom{n}{k}.
\end{equation}
Note that for each pair $(S, c)$ where $S$ is a set of $k-1$ vertices and $c$ is a color, we have exactly one of the following:
\begin{enumerate}[i)]
    \item there is no $c$-colored edge containing $S$, 
    \item there is some $c$-colored edge $e \supseteq S$, but there is no other $c$-colored edge intersecting $e$ in $k-1$ vertices, or
    \item \label{case:2comp} there are some $c$-colored edges $e_1,e_2$ with $e_1 \supseteq S$ and $|e_1 \cap e_2|=k-1$.
\end{enumerate}
Note that in case \ref{case:2comp}), there cannot be a third $c$-colored edge $e_3$ with $|e_3 \cap (e_1 \cup e_2) |\ge k-1$. The total number of pairs $(S, c)$ of a $(k-1)$-set of vertices and a color is
\begin{equation}\label{eqn:paircount}
x_0 + kx_1 + (2k-1)x_2 = \binom{n}{k-1} |C|.
\end{equation}
Now note that for every pair $(e_1 \cup e_2, c)$ counted by $x_2$,  there are $k-1$ pairs $(e', c')$ counted by $x_1$. Indeed, there are $k-1$ edges $e'\neq e_1, e_2$ contained in $e_1 \cup e_2$. Each such edge $e'$ has some color $c'$ and there cannot be any other $c'$-colored edge $e''$ intersecting $e'$ in $k-1$ vertices, since otherwise $e_1 \cup e_2 \cup e''$ would be have two repeated colors. Furthermore such an edge $e'$ cannot be contained in the union of two other edges, say $e_1', e_2'$ of the same color since then $e_1, e_2, e_1', e_2'$ would span $k+2$ vertices and have two repeated colors. Thus
\begin{equation}\label{eqn:ybound}
   x_1  \ge (k-1)x_2. 
\end{equation}
But now, as explained below, we have
\begin{align*}
    \binom{n}{k-1} |C| & = x_0 + kx_1 + (2k-1)x_2\\
    & \ge \left(k-\frac{1}{k+1}\right)(x_1 + 2x_2) + \frac{1}{k+1} \bigg( x_1 - (k-1)x_2 \bigg)\\
    & \ge \left(k-\frac{1}{k+1}\right) \binom{n}{k}.
\end{align*}
Indeed, the first line is just \eqref{eqn:paircount}. The second line follows since $x_0 \ge 0$ and since the coefficients of $x_1, x_2$ are the same as in the line above. The last line follows from  \eqref{eqn:edgecount} and \eqref{eqn:ybound}. Solving for $|C|$ above yields 
\[
|C| \ge \frac{\left(k-\frac{1}{k+1}\right) \binom{n}{k}}{\binom{n}{k-1}}=  \frac{k^2+k-1}{k^2+k}  n + O(1).
\]

\subsection{Upper bound}

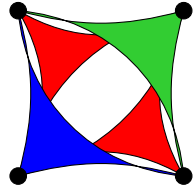
\begin{figure}[h]
        \centering
        \begin{tikzpicture}[scale=2.2]
         \draw[fill = red, opacity = 0.5] (0,0) to[bend left = 30] (1,1) to[bend left = 30] (0,1) to [bend left = 30] cycle;

\draw[fill = red, opacity = 0.5] (0,0) to[bend left = 30] (1,0) to[bend left = 30] (1,1) to [bend left = 30] cycle;

\draw[fill = green, opacity = 0.5] (0,1) to[bend right = 15] (1,1) to[bend right = 15] (1,0) to [bend right = 40] cycle;

\draw[fill = blue, opacity = 0.5] (0,0) to[bend left = 15] (1,0) to[bend left = 40] (0,1) to [bend left = 15 ] cycle;

\filldraw (0,0) circle (0.05 cm);

\filldraw (1,0) circle (0.05 cm);

\filldraw (0,1) circle (0.05 cm);

\filldraw (1,1) circle (0.05 cm);

        \end{tikzpicture}
        \caption{A tile  in the case $k=3$.}
        \label{fig:tile}
    \end{figure}

As before, we color in two stages. In the first stage, we partially color $K_n^k$ with $k$-colored copies of $K_{k+1}^k$ in which every $k$-tuple of vertices receives a distinct color, except for two which receive the same color, creating one repetition. We will refer to such colored copies of $K_{k+1}^k$ as \emph{tiles} in this coloring (as shown in Figure~\ref{fig:tile} for $k=3$). We obtain this partial coloring using the conflict-free hypergraph matching method.  

A path of length 2 in a $k$-uniform hypergraph is a pair of edges $e,e'$ that have $k-1$ vertices in common. 

\begin{theorem}\label{thm:cliquesproperties}
There exists $\delta>0$ such that for all sufficiently large $n$ in terms of $\delta$, there is an edge-coloring of a subgraph $F\subset K_n^k$ with at most $\frac{k^2+k-1}{k^2+k}n+n^{1-\delta}$ colors and the following properties:
\begin{enumerate}
    \item Every color class is a union of edges and paths of length 2, where any two such components intersect in at most $k-2$ vertices. \label{prop:clique1}
    \item In any tile $K_{k+1}^k$ in $F$ in which edges $e$ and $e'$ receive the same color $i_1$ and all other edges receive distinct colors $i_2,\ldots,i_k$, the $(k-1)$-set of vertices $e\cap e'$ is not contained in any edge of color $i_j$ for $2\leq j\leq k$.  \label{prop:clique2}
    \item The coloring on $F$ is a $\left(K_{k+2}^k,\binom{k+2}{k}-1\right)$-coloring. \label{prop:clique3}
    \item Any $(k-1)$-set of vertices belongs together to at most $n^{1-\delta}$ edges  in $L=K_n^k\setminus E(F)$.  \label{prop:clique4}
    \item Let $e\in E(K_n^k)$ be $e=\{x,y\}\cup S$ for some $S\in\binom{n}{k-2}$. Then the number of edges $e'\in E(L)$ with $e'=\{x',y'\}\cup S$ such that $\{x,x'\}\cup S$ and $\{y,y'\}\cup S$  receive the same color is at most $n^{1-\delta}$. \label{prop:clique5}
\end{enumerate}
\end{theorem}

\begin{proof}
Let $\delta>0$ be sufficiently small, and let $n$ be sufficiently large in terms of $\delta$. Let $\rho=n^{-\delta}$, and set the number of colors $c=(1+\rho)(\frac{k^2+k-1}{k^2+k})n$. We will define a random set $V$ as follows. Let $V'=\bigcup_{i\in [c]}V_i'$, where $V_1',\ldots,V_c'$ are copies of $\binom{[n]}{k-1}$. Delete each element in $V'$ independently with probability $p=\rho/(1+\rho)$, and let $V_i$ denote the remaining elements in $V_i'$. Let $V=\bigcup_{i\in [c]} V_i$. For each color $i\in [c]$ and $S\in\binom{[n]}{k-1}$, denote by $S_i$ the copy of $S$ in $V_i$, if it exists.

We now construct a hypergraph $\HH$ with vertex set $E(K_n^k)\cup V$.  
For each set $B\in\binom{[n]}{k+1}$, a choice of $k$ colors $i_1,i_2,\ldots,i_{k+1}\in[c]$ where $i_1=i_2$, and an ordering $B_1,\ldots, B_{k+1}$ of the $k$-tuples of vertices in $\binom{B}{k}$, we add to $\HH$ the edge
\[\{B_1,\ldots,B_{k+1}\}\cup \left\{S_{i_j}:S\in\binom{B_j}{k-1}\right\}.\] 
This corresponds to a tile, which is a copy of the clique $K_{k+1}^k$ on vertex set $B$ in which every edge $B_j$ receives the color $i_j$, and in particular, $B_1$ and $B_2$ receive the same color. 
Thus, any matching in $\HH$ corresponds to a collection of edge-disjoint tiles which share no $(k-1)$-tuple in the same color. This yields a coloring of a subhypergraph of $K_n^k$ satisfying properties \ref{prop:clique1} and \ref{prop:clique2}. In order to guarantee properties \ref{prop:clique3}-\ref{prop:clique5}, we will later define an appropriate conflict system and trackable test functions for $\HH$. 

First, however, we check that the degree conditions in Theorem~\ref{thm:blackbox} hold for $\HH$. 
Note that $1-p=(1+\rho)^{-1}$. 
For a vertex in $\HH$ of the form $B\in \binom{[n]}{k}$, 
\begin{align*}
\E[\deg_{\HH}(B)]&=(n-k)\cdot c^{\underline{k}} \cdot \left(k+\binom{k}{2}\right)\cdot (1-p)^{2\binom{k+1}{k-1}-1}p^{(k-2)\binom{k+1}{k-1}+1}\\
&=\frac{k^2+k-1}{2}n^{k+1}(1-p)^{2\binom{k+1}{k-1}-2}p^{(k-2)\binom{k+1}{k-1}+1}\pm O(n^{k}).
\end{align*}
For a vertex in $\HH$ of the form $S_i\in V_i'$, 
\begin{align*}
\E[\deg_{\HH}(S_i):S_i\in V_i]&=(c-1)^{\underline{k-1}}\left(\binom{n-k+1}{2}+(n-k+1)(n-k)+2(n-k+1)(n-k)\right)\\
&=\frac{k^2+k-1}{2}n^{k+1}(1-p)^{2\binom{k+1}{k-1}-2}p^{(k-2)\binom{k+1}{k-1}+1}\pm O(n^{k}).
\end{align*}
Applying Theorem~\ref{thm:McDiarmid} shows that with high probability, 
\begin{equation*}
\begin{split}
        \frac{k^2+k-1}{2}n^{k+1}  (1-p)^{2\binom{k+1}{k-1}-2}& p^{(k-2)\binom{k+1}{k-1}+1}-O(n^{k+\frac{2}{3}}) \\
        &\leq \min_{u\in V(\HH)} \deg_{\HH}(u) \\& \leq \max_{u\in V(\HH)} \deg_{\HH}(u) 
        \\&\leq \frac{k^2+k-1}{2}n^{k+1}(1-p)^{2\binom{k+1}{k-1}-2}p^{(k-2)\binom{k+1}{k-1}+1}+ O(n^{k+\frac{2}{3}}).
\end{split}
\end{equation*}

To see this, fix a vertex in $\HH$, say $B\in \binom{[n]}{k}$. For $S_i\in V'$, define $b_{S_i}=nc^{k-1}$ if $S\subset B$, $b_{S_i}=c^{k-1}$ if $|S\cap B|=k-2$, and $b_{S_i}=0$ otherwise. Thus, $\sum_{S_i\in V'}b_{S_i}^2=O(n^{2k+1})$, and hence we have concentration on an interval of length $O(n^{k+\frac{1}{2}})$ by Theorem~\ref{thm:McDiarmid}. Now instead fix a vertex $S_j'\in V_j$, and consider $S_i\in V'$ with $S_i\neq S_j'$. Set $b_{S_i}=n^k$ if $S=S'$ or if $|S\cap S'|=k-2$ and $i=j$, set $b_{S_i}=n^{k-1}$ if $|S\cap S'|=k-2$ and $i\neq j$ or if $|S\cap S'|=k-3$ and $i=j$, set $b_{S_i}=c^{k-2}$ if $|S\cap S'|=k-3$ and $i\neq j$, and set $b_{S_i}=0$ otherwise. Then $\sum_{S_i\in V'}b_{S_i}^2=O(n^{2k+1})$, as before, and the same conclusion holds. Thus, $\HH$ is almost $d$-regular with $d=\frac{k^2+k-1}{2}n^{k+1-\delta\left((k-2)\binom{k+1}{k-1}+1\right)}$. 

Moreover, by a similar application of Theorem~\ref{thm:McDiarmid}, we have that with high probability, $\Delta_2(\HH)=O(n^k)$. From now on, we fix a deterministic hypergraph that satisfies these properties and refer to this as $\HH$. 

Now we define a conflict system $\CC$ for $\HH$. Let $V(\CC)=E(\HH)$, and let the edges of $\CC$ correspond to copies of $K_{k+2}^{k}$ with at most $\binom{k+2}{k}-2$ colors. Note that every conflict in $\CC$ has size 4; indeed, first notice that by the definition of the tiles in $\HH$ and the fact that we get a matching, there cannot be a copy of  $K_{k+2}^k$ with 3 $k$-tuples receiving the same color $i$. Therefore, if we have a $K_{k+2}^k$ with at most $\binom{k+2}{k}-2$ colors, then there are two distinct colors $i$ and $j$ which are repeated, each coming from a separate tile, creating a conflict of size 4. This shows that conflicts of size 4 with 2 repeated colors are enough to forbid any bad $K_{k+2}^k$. 


Condition (C2) is satisfied since $\Delta(\CC)\leq O(n^2 \cdot d\cdot \frac{d}{n}\cdot \frac{d}{n})=O(d^3)$. Indeed, given a tile in $V(\CC)$, there are $O(n^2)$ ways to pick the other two graph vertices in the conflict, $O(d)$ ways to pick a tile that contributes a different color to the conflict, and $O(d^2/n^2)$ ways to pick the remaining two tiles in the conflict. 

Similarly, (C3) is satisfied since $\Delta_2(\CC)\leq O(d^{2-\eps})$ and $\Delta_3(\CC)\leq O(d^{1-\eps})$ for all $\eps\in(0,\frac{1}{k+2})$. Indeed, a pair of tiles contained in a conflict either contribute different colors or the same color, and in both cases, there are $O(d^2/n)$ ways to complete the conflict. Similarly, given three tiles, there are $O(d/n)$ ways to complete a conflict containing all three. Thus, $\CC$ is $(d,\ell,\eps)$-bounded for $\eps\in (0,\frac{1}{k+2})$. 

Applying Theorem~\ref{thm:blackbox} with $\HH$ and $\CC$ would yield a conflict-free matching $M$ in $\HH$ which would correspond to a coloring of a subgraph $F$ of $K_n^k$ satisfying properties \ref{prop:clique1}-\ref{prop:clique3} of Theorem~\ref{thm:cliquesproperties}. 
In order to obtain a coloring which also satisfies properties \ref{prop:clique4} and \ref{prop:clique5} of Theorem~\ref{thm:cliquesproperties}, we now introduce appropriate test functions. 

First, we will define a 1-uniform test function for each $(k-1)$-tuple $X$ of vertices in $V(K_n^k)$. Let $S_X$ be the set of edges containing $X$ in $K_n^k$. Let $w_X:E(\HH)\rightarrow\{0,2\}$ be the weight function which assigns to each edge $e\in E(\HH)$ the size of its intersection with $S_X$. Then we have  
\[w_X(\HH)=\sum_{e\in S_X}\deg_{\HH}(e)=nd-O(n^{k-2})\geq d^{1+\eps},\]
satisfying condition (W1). In addition, conditions (W2)-(W4) are trivially satisfied since $w_X$ is 1-uniform, so $w_X$ is a $(d,\eps,\CC)$-trackable  test function for all $\eps\in(0,\frac{1}{k+2})$ and $v\in V(K_n)$. 

Now fix $\delta<\eps^3\log_n(d)$. Applying Theorem~\ref{thm:blackbox} with these test functions yields a $\CC$-free matching $M$ in $\HH$ for which 
\[w_X(M)>(1-d^{-\eps^3})d^{-1}w_X(\HH)>(1-n^{-\delta})n\]
for each $(k-1)$-tuple $X$ of vertices in $V(K_n^k)$. 
Consequently, any such $(k-1)$-tuple of vertices has codegree at most $n-(1-n^{-\delta})n=n^{1-\delta}$ in the uncolored subgraph $L=K_n^k-E(F)$, proving property \ref{prop:clique4} of Theorem~\ref{thm:cliquesproperties}.

In order to also obtain property \ref{prop:clique5} of Theorem~\ref{thm:cliquesproperties}, we will define two more types of test functions for each edge $e\in E(L)$. This will allow us to bound the number of bad colorings of $K_{k+2}^k$ which could contain the edge $e$ and two edges in $F$. For the first test function, let $\PP_{e}$ be the set containing any matching of two tiles $B,B'$ for which the edges $f\in B$ and $f'\in B'$ receive the same color $i\in[c]$, and $|f\cap f'|=k-2$, and $|e\cup f\cup f'|=k+2$. Note that there are three possibilities for how $e$, $f$, and $f'$ can be arranged in such a matching. 

Let $w_e$ be the indicator weight function for this set $\PP_e$. This is a 2-uniform test function for $\HH$, and we will show that $w_e$ is $(d,\eps,\CC)$-trackable for each $\eps\in(0,\frac{1}{k+2})$. To see that condition (W1) holds for $w_e$, note that there are $\binom{n}{2}$ ways to pick the vertices $x'$ and $y'$ which are in $(f\cup f')\backslash e$, $(\binom{k}{k-1}\binom{k-1}{k-2}+\binom{k}{k-2}\binom{k-2}{k-4}+2\binom{k}{k-1}\binom{k-1}{k-3})$ ways to complete pick $e$ and $f$ from the set of $k+2$ vertices $e\cup\{x,y\}$, $cp^{2k}$ ways to pick a color for $f$ and $f'$, and $d/c$ ways each to pick the tiles $B$ and $B'$ containing $f$ and $f'$, respectively:
\[w_e(\HH)=|\PP_e|\geq \binom{n}{2}\cdot \left(k(k-1)+\binom{k}{k-2}\binom{k-2}{k-4}+2(k-1)\binom{k-1}{k-3}\right)\cdot cp^{2k}\cdot \left(\frac{d}{cp^k}\right)^2\pm O(n^{2k+2})>d^{2+\varepsilon}.\]

Now fix a tile $B\in E(\HH)$. There are at most $O\left(n\cdot \frac{d}{cp^k}\right)<\frac{w_e(\HH)}{d^{1+\eps}}$ pairs in $\PP_e$ which contain $B$ for all $\eps\in(0,\frac{1}{k+2})$, so condition (W2) is satisfied. 

To verify condition (W3), fix two tiles $B,B'\in E(\HH)$ which form a pair in $\PP_e$. To count the triples of tiles which form a conflict with $B$, we must choose at most 5 other vertices outside $B\cup B'$ and $1+3(k-1)=3k-2$ colors. This gives $O(n^{3k+3})$ choices for these triples, however, since the conflicts must also be completed using $B'$ instead of $B$, we save either two choices of vertices or a choice of vertex and a choice of color; hence, there are at most $O(n^{3k+1})<d^{3-\eps}$ ways to complete a conflict with $B$ or $B'$, so condition (W3) holds. 

Finally, condition (W4) holds vacuously, so $w_{e}$ is $(d,\eps,\CC)$-trackable for each $e\in E(L)$. 

Next, for each $e\in E(L)$, we define a second test function $w'_e$ which is an indicator weight function for the set $\TT_e$ of triples of tiles which extend the pairs in $\PP_e$ as follows. Given a pair $B,B'$ of tiles in $\PP_e$ containing edges $f$ and $f'$ of color $i$, respectively, we add the triple $\{B,B',B''\}$ to $\TT_e$ whenever the tile $B''$ is edge-disjoint from $B$ and $B'$ and contains an edge $e'\notin\{e,f,f'\}$ on some $k$ of the $k+2$ vertices in $e\cup f \cup f'$ in some color $j\in[c]\backslash \{i\}$.

This is a 3-uniform test function for $\HH$, and  we will show that it is $(d,\eps,\CC)$-trackable for $\eps\in(0,\frac{1}{k+2})$. Note that 
\[w'_e(\HH)=|\TT_e|=\left(\binom{k+2}{k}-3\right)(d\pm O(n^{k}))|\PP_e|>d^{3+\eps},\]
so condition (W1) is met for $w_e'$.

To check condition (W2), first fix a tile $B\in E(\HH)$. There are at most $O\left(n\cdot \frac{d}{c}\cdot d\right)<w'_e(\HH)/d^{1+\eps}$ triples in $\TT_e$ which contain $B$ for all $\eps\in(0,\frac{1}{k+2})$. Now fix a pair of tiles $B$ and $B'$. There are at most $O(\left(d\right)<w'_e(\HH)/d^{2+\eps}$ triples in $\TT_e$ containing $B$ and $B'$ for all $\eps\in(0,\frac{1}{k+2})$.

In addition, to see that condition (W3) is satisfied for $w'_e$, fix a pair of tiles $B,B'$ which appear in a triple in $\TT_e$ together. Then, as for $w_e$, there are at most $O(n^{3k+1})<O(d^{3-\eps})$ triples of tiles which can form a conflict with $B$ and with $B'$. 

Finally, as before, condition (W4) is vacuously true, and $w'_e$ is $(d,\eps,\CC)$-trackable for all $\eps\in(0\frac{1}{k-1})$.

We can now apply Theorem~\ref{thm:blackbox} with all of the test functions of type $w_X$, $w_e$, and $w'_e$ to obtain a $\CC$-free matching $M$ in $\HH$ for which property \ref{prop:clique4} of Theorem~\ref{thm:cyclesproperties} holds and we have 
\[w_e(M)=\left|\binom{M}{2}\cap \PP_e\right|\leq (1+d^{-\eps^3})d^{-2}w_e(\HH)\leq (1+n^{-\delta})n^{-2k-2}|\PP_e|\] 
and 
\[w'_e(M)=\left|\binom{M}{3}\cap \TT_e\right|\geq (1-d^{-\eps^3})d^{-3}w'_e(\HH)\geq (1-n^{-\delta})n^{-3k-3}\left(\binom{k+2}{k}-3\right)|\PP_e|\]
for all $e\in E(L)$.
Thus, the number of edges $e'\in E(L)$ as described in property \ref{prop:clique5} of Theorem~\ref{thm:cyclesproperties} is at most the difference 
\[\left(\binom{k+2}{k}-3\right) w_e(M)-w'_e(M)\leq n^{1-\delta},\]
which completes the proof of Theorem~\ref{thm:cliquesproperties}.
\end{proof}

We can now apply the Local Lemma to color the edges of $L$ using a new set of $n^{o(1)}$ colors without creating any bad copies of $K_{k+2}^k$, completing the proof of Theorem \ref{thm:hypergraphs}. Since this argument is very similar to those appearing in Section \ref{sec:cycles} as well as in \cite{BHZ,JM}, we omit the details.

\section{Paths in complete graphs}\label{sec:paths}

Here we prove Theorem \ref{thm:paths}, starting with the proof of (i), that is, $f(K_n, P_6, 4)=\frac 14 n^2 + O(n)$. 
First, to prove the upper bound on $f(K_n, P_6, 4)$, we construct a coloring as follows. We start by packing edge-disjoint copies of $K_4$ into $K_n$. With the right divisibility conditions we can get a perfect packing, but in any case we can pack at least $\frac{1}{12}n^2 -O(n)$ copies of $K_4$ (i.e. we cover all but $O(n)$ edges of $K_n$). Now we will use 3 colors to properly color each copy of $K_4$ in the packing, such that each $K_4$ has 3 colors that are not used anywhere else. For the $O(n)$ edges of $K_n$ that are not covered by the packing, we give each one its own color that is used nowhere else. Altogether this uses $\frac 14 n^2 + O(n)$ colors. 

To see that this coloring is a $(P_6, 4)$-coloring, suppose we have a copy of $P_6$ with 2 repeats. This can only happen if the $P_6$ has 2 edges of some color $c_1$ and 2 edges of some other color $c_2$. The $c_1$-colored edges and $c_2$-colored edges  cannot belong to the same $K_4$ in the packing, or else these 4 edges would make a cycle. But 2 distinct copies of $K_4$ in the packing would induce at least 7 vertices, so we have a contradiction. Thus the coloring is a $(P_6, 4)$-coloring, showing $f(K_n, P_6, 4)\leq \frac{1}{4}n^2+O(n)$.

We now prove the lower bound. Suppose we have a $(P_6, 4)$-coloring of $K_n$. We claim that we can remove at most 3 vertices to obtain a proper coloring of $K_{n'}$, where $n-3 \le n' \le n$. Indeed, let $\GG$ be a 3-uniform hypergraph with vertex set $V(K_n)$ and an edge $e_1 \cup e_2$ whenever $e_1, e_2 \in E(K_n)$ are adjacent edges of the same color. Then $\GG$ cannot have a matching with 2 edges, since those 6 vertices would contain a $P_6$ with 2 repeats. Thus, we can remove all 3 vertices from any edge of $\GG$ to destroy all edges of $\GG$. This proves our claim. 

    So we have a $(P_6, 4)$-coloring of $K_{n'}$ which is also a proper coloring. Now note that if any color is used 3 times then we easily find a $P_6$ with 2 repeats. Thus each color is used at most twice and the number of colors is at least $\binom{n'}{2}/2 = \frac 14 n^2 + O(n)$. Thus $f(K_n, P_6, 4) \ge \frac 14 n^2 + O(n)$.

Now we prove part (ii) of Theorem \ref{thm:paths}, that is,  $f(K_n, P_7, 5)= \frac 12 n^2 +O(n)$. Of course, the upper bound is trivial. 
For the lower bound, suppose we have a $(P_7, 5)$-coloring of $K_n$. This coloring must also be a $(P_6, 4)$-coloring, and so as in the last proof we can remove at most 3 vertices to obtain a proper coloring. So assume now that we have a proper $(P_7, 5)$-coloring of $K_{n'}$ with $n-3 \le n' \le n$. Now observe that each color is used at most twice. Furthermore if $e_1, e_2$ have the same color $c$ and $e_3, e_4$ have the same color $c'$, then $e_1 \cup e_2$ and $e_3 \cup e_4$ are disjoint. Thus, only $O(n)$ colors can be used twice. This completes the proof of our estimate for $f(K_n, P_7, 5)$.

\begin{figure}[h]
\centering
\begin{tikzpicture}[scale=1]
\foreach \i in {0,1,2}{
	\node[inner sep=.08cm,fill,circle] (a\i) at (120+\i*120:1.5) {};
	\node[inner sep=.08cm,fill,circle] (b\i) at (60+\i*120:1.5) {};
    \draw[ultra thick, red] (a\i)--(b\i);
    }
\node at ([shift={(90:.35)}]a0) {$a_1$};
\node at ([shift={(90:.35)}]b0) {$b_1$};
\node at ([shift={(270:.35)}]a1) {$a_3$};
\node at ([shift={(180:.35)}]b1) {$b_3$};
\node at ([shift={(0:.35)}]a2) {$a_2$};
\node at ([shift={(270:.35)}]b2) {$b_2$};
\draw[ultra thick, blue] (a0)--(b2) (b0)--(a2); 
\draw[ultra thick, green!50!teal] (a0)--(a2) (b0)--(b2);
\end{tikzpicture}
\caption{A portion of the coloring for each $K_6$ in the packing for the upper bound on $f_{\text{proper}}(K_n, P_8, 5)$}
\label{fig:fP85}
\end{figure}
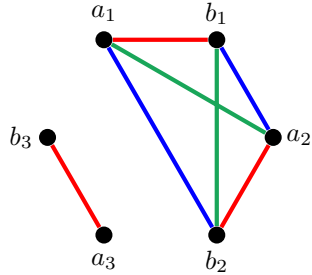

Finally we prove (iii), that is,  $f_{\text{proper}}(K_n, P_8, 5) = \frac{7}{30}n^2 +O(n)$. First we prove the upper bound. To obtain a coloring, we start with a packing of $\binom n2 / \binom 62 -O(n) = \frac 1{30} n^2 - O(n)$ edge-disjoint copies of $K_6$. On each copy of $K_6$ in the packing we will use 7 colors which are not used anywhere else besides on the edges of this $K_6$. Let the vertices of $K_6$ be $\{a_1, b_1, a_2, b_2, a_3, b_3\}$. Color the edges $a_ib_i$ with color 1 for all $i\in[3]$. For $i,j\in[3]$, with $i\neq j$, we properly color the 4 edges between $\{a_i,b_i\}$ and $\{a_j, b_j\}$ using 2 new colors. See Figure \ref{fig:fP85} where a portion of the coloring is shown. For the $O(n)$ edges of $K_n$ that are not covered by the packing, we give each one its own color that is used nowhere else. Altogether this uses $\frac{7}{30} n^2 + O(n)$ colors.

To prove that this is a $(P_8, 5)$-coloring, suppose we have a $P_8$ with $3$ repeats. In the arguments that follow, let $A_i$ represent the set of vertices incident to $c_i$-colored edges. 
First suppose color $c_1$ appears 3 times and color $c_2$ appears 2 times (on our path $P_8$). So $|A_1|=6$ and $|A_2|=4$. The $c_1$-colored edges and the $c_2$-colored edges cannot belong to the same $K_6$ in the packing or else they would form a cycle. So $A_1$ and $A_2$ must belong to different copies of $K_6$ in the packing and hence $|A_1\cap A_2| \le 1$.
Thus we have $|A_1\cup A_2| = |A_1| + |A_2| - |A_1\cap A_2| \ge 4+6-1 = 9>8$, a contradiction. 
Now suppose there exists a $P_8$ in which colors $c_1, c_2$ and $c_3$ each appear twice, so $|A_1| = |A_2| = |A_3|=4$. These edges cannot all appear on the same $K_6$ in the packing else we would have 6 edges on 6 vertices which would result in a cycle.  If colors $c_1$ and $c_2$ appear on the same $K_6$ in the packing then they must span all 6 vertices and so we would have $|(A_1\cup A_2) \cup A_3| = 6+4 - 1 = 9>8$. Finally, if colors $c_1, c_2$ and $c_3$ all appear on distinct copies of $K_6$ in the packing, then we have $|A_1\cup A_2 \cup A_3| \ge 3\cdot4 - \binom32 \cdot 1 = 9>8$.

To prove the lower bound, suppose we have a proper $(P_8,5)$-coloring of $K_n$. We note that if any color is used 4 times, then we easily find a $P_8$ with 3 repeats. For $i\in[3]$, we let $x_i$ represent the number of colors which appear $i$ times. Then we have $x_1 + 2x_2 + 3x_3 = \binom{n}{2}$. Suppose color $c_1$ is used 3 times and let $A_1$ represent the 6 vertices spanned by these edges. Any color used among the vertices of $A_1$ cannot be incident to any vertex outside of $A_1$ otherwise we can find a $P_8$ with 3 repeats. Also, any color other than $c_1$ appearing among $A_1$ can appear at most twice. Thus, for every color appearing 3 times, the copy of $K_6$ it spans has 12 other edges whose colors must all appear either once or twice (and if they appear twice, both instances are within this same $K_6$). Thus we have the inequality $x_1 + 2x_2 \ge 12x_3$. Therefore, the number of colors used is
    \begin{align*}
        x_1 + x_2 + x_3 &= \frac{7}{15}(x_1+2x_2+3x_3) + \frac{1}{30}(x_1 + 2x_2 - 12x_3) + \frac12 x_1 \\
        &\ge  \frac{7}{15}(x_1+2x_2+3x_3) \\
        &= \frac{7}{15}\binom{n}{2}.
    \end{align*}

We were unable to asymptotically determine $f(K_n, P_8, 5)$. The problem we have is that we cannot seem to emulate the first ``cleaning'' step of the previous two proofs, i.e. removing a few vertices to obtain a proper coloring. We strongly suspect that $f(K_n, P_8, 5) = \frac{7}{30}n^2 +O(n)$ as well.

\section{Conclusion}\label{sec:conclusion}

Using similar methods, we can obtain new bounds on the generalized Ramsey number of a tight cycle in a complete $k$-uniform hypergraph. Let $C_{\ell}^k$ denote a tight $k$-uniform cycle on $\ell$ vertices, where consecutive edges in the cycle intersect in $k-1$ vertices. Then we have
\[\frac{2}{k \ell + \ell - 1}n\leq f(K_n^k,C_{\ell}^k,k+1)\leq \frac{1}{\ell-k}n+o(n).\]
Indeed, the upper bound comes from using the conflict-free hypergraph matching method to find a large packing of monochromatic copies of $K_{\ell-1}^k$ where any two copies of the same color overlap in at most $k-2$ vertices (and where we avoid some additional conflicts).  
To prove the lower bound, observe that no color class in a $(C_{\ell}^k,k+1)$-coloring of $K_n^k$ can contain a monochromatic path $P^k_{\ell}$. Now we apply a result of~\cite{FJKMV}  giving an upper bound on the extremal number of this path to get a lower bound on the number of colors in any $(C_{\ell}^k,k+1)$-coloring of $K_n^k$. Indeed, since $\ex(K_n^k,P_{\ell}^k)\leq \frac{k \ell + \ell - 1}{2k}\binom{n}{k-1}$, we have
$f\left(K_n^k,C_{\ell}^k,k+1\right)\geq \binom{n}{k}\cdot \frac{2k}{(k \ell + \ell -1)\binom{n}{k-1}}$. 

\emph{Remark:} While the writing of this manuscript was nearing completion, we learned that results similar to our Theorems \ref{thm:cycles} and \ref{thm:bipartitecycles}
 were obtained independently by Andrew Lane and Natasha Morrison \cite{LM}. 

\bibliography{bibfile}
\bibliographystyle{abbrv}

\end{document}